\let\bsbf=\bf
\let\bsit=\it
\let\bs=\rm
\documentstyle[12pt]{article}

\def\triangulo{\hbox{.\kern2pt.\kern-5pt\raise 4pt\hbox{.}}\kern 5pt}
\font\gordas = msbm10 at 12pt

\def\bbb#1{\hbox {{\gordas #1}}}

\def\a{{\bbb A}}

\def\ce{{\bbb C}}
\def\ache{{\bbb H}}
\def\o{{\bbb O}}
\def\erre{{\bbb R}}

\begin{document}
\bs
\begin{center}
{\large\bsbf The zero divisors of the Cayley--Dickson}\\
{\large\bsbf algebras over the real numbers} \\[.5cm]
Guillermo Moreno
\end{center}
\footnote{Classification 1991 AMS 17A99}
\footnote{Accepted for publication at Bol. Soc. Mat. Mex.}
\vglue1cm
\noindent
{\bsbf Abstract. } In this paper we describe algebraically the zero divisors of the Cayley - Dickson algebras
$\a_{n} = \erre^{2^n}$ for $n\ge 4$ over the real numbers.
\vglue1cm
\noindent
{\bsbf Introduction.} The Cayley--Dickson algebra $\a_{n}$ 
over $\erre$ 
is an algebra structure on $\erre^{2^n}$ given inductively by the 
formulae:
\par 
Let $x=(x_{1},x_{2})$ and $y=(y_{1},y_{2})$ in 
$\erre^{2^n}=\erre^{2^{n-1}}\times\erre^{2^{n-1}}$ then
$$xy=(x_{1}y_{1}-\overline{y}_{2}x_{2},y_{2}x_{1}+x_{2}\overline{y}_{1})$$
where
$$\overline{x}=(\overline{x}_{1},-x_{2}).$$
Therefore $\a_{0}=\erre, \a_{1}=\ce$ complex numbers, $\a_{2}=\ache$ the
Hamilton quaternions, $\a_{3}=\o$ the Cayley octonians, etc. This four algebras:
$\erre, \ce, \ache$ and $\o$ are known as the classical Cayley--Dickson algebras
and their main distinctive feature of these is:
\vglue.5cm
\noindent
{\bsbf Hurwitz Theorem:} Let ``$\parallel\quad\parallel$'' denote the euclidean 
norm in $\erre^{2^n}$. Then
\par 
$\parallel xy\parallel=\parallel x\parallel \parallel y\parallel \forall x$
and $y$ in $\a_{n}$ if and only if $n=0,1,2,3$. (See [5] and [7]).
\par 
That is, for $n\ge 4$ this norm--preserving formula is not 
true in general and this 
opens the possibility of the existence of \underline{zero divisors} in $\a_{n}$ for
$n\ge 4$ i.e. $x$ and $y$ non--zero elements in $\a_{n}$ such that
$xy=0$ e.g. let $x=e_{1}+e_{10}$ and $y=e_{15}-e_{4}$ in $\erre^{16}=\a_{4}$ where
$e_{0},e_{1},e_{2},\ldots,e_{15}$ is the canonical basis in $\erre^{16}$.
\par 
In this paper we study  the zero divisors in 
$\a_{n}$ for $n\ge 4$.
\par 
From the algebraic side these algebras are non--commutative for $n\ge 2$ and
non--associative for $n\ge 3$. Moreover $\a_{3}$ is \underline{alternative}:
\par 
$x^{2}y=x(xy)$ and $xy^{2}=(xy)y\,\,$ for all $x$ and $y$ in $\a_{3}$

and $\a_{n}$ for $n\ge 4$ is \underline{flexible}:
\par 
$x(yx)=(xy)x\,\,$ for all $x$ and $y$ in $\a_{n}$. Clearly
$$\hbox{\bs Associative}\quad\Rightarrow\quad
\hbox{\bs Alternative}\quad\Rightarrow\quad\hbox{\bs Flexible}$$
and the backwards implications are not true. Introducing the ``associator 
symbol''
$$(a,b,c)=(ab)c-a(bc)$$
we have that alternativity property in $\a_{3}$ is equivalent to the 
alternativity of the associator symbol in $\a_{3}$ via polarization; flexibility 
means that $(x,y,x)=0$ for all $x$ and $y$ in $\a_{n}$.
\par 
We will see that the absence of the norm--preserving property, the 
non--alternativity and the presence of zero divisors in $\a_{n}$ for $n\ge 4$ are very
strongly related and in fact they determine one to each other. 
\par 
Usually the zero 
divisors are present in different algebraic contexts
\par 
\begin{itemize}                                               
\item[1)]the trivial one: in the direct sum of algebras i.e. if 
${\cal A}$ and ${\cal B}$
are algebras over $\erre.\quad {\cal A} \oplus {\cal B}$ has a lot
of zero divisors e.g. the ones of the form $(a,0)$ and $(0,b)$.
\item[2)]In an associative (alternative) context as in 
$M_{n}(\erre)$ the $n\times n$
matrices over $\erre$ where a god--given invariant called  
determinant says that \\
$0\neq A\in M_{n}(\erre)$ is a zero divisor if and only if $det A=0$.
\item[3)]In a non--associative (non--alternative) context which becomes rather 
difficult, this is our case.
\end{itemize}
\vglue.5cm
This paper has two chapters, in chapter one we give a general description of 
the zero divisors for $\a_{n}\quad (n\ge 4)$ studying the linear transformations
$R_{a},L_{a}:\a_{n}\rightarrow \a_{n}$ right and left multiplication respectively for 
$0\neq a \in \a_{n}$ fixed and  we prove that each $a\neq 0$ (doubly pure)
induce a direct sum 
decomposition of $\a_{n}$ where one summand is a copy of $\ache$ the 
quaternions, 
other summand is the elements of $\a_{n}$ which ``alternate'' with $a$ and one
third
summand is the annihilator of $\a$ ``$Ker L_{a}=Ker R_{a}$'' (see Theorem I.15).
\par 
This also implies that $dim Ker L_{a}\equiv 0$ mod 4 and $dim Ker L_{a}\le 2^{n}-4$.
\par 
In Chapter II we construct some type of zero divisors (called special ones), in 
$\a_{4}$ are all of them, which are pairs of alternatives elements, 
we 
describe all of them completely for $n\ge 4$ 
and relate these zero divisors with the 
multiplicative monomorphisms of $\a_{3}$ to $\a_{n}\quad n\ge 3$.
\par 
In a second paper we will describe the topology of the subset of 
$\erre^{2^{n+1}}$%
$$
ZD^{n}_{1}=\{(x,y)\in \erre^{2^n}\times \erre^{2^n}|xy=0\,\,\hbox{\bs and}\,\,
|x|=|y|=1\}$$
and relate this with suitable Stiefel manifolds. (Second paper is under review, the preprint is available under
request).
\par 
We work also on a third paper on the subject: ``Applications'' where we use 
the algebraic and topological parts to construct bilinear normed maps, 
construct generalized Hopf maps, etc.
\par 
We thank to Fred Cohen, Sam Gitler and K.Y. Lam for many illuminating
conversations and to Isidoro Gitler and Jos\'e Mart{\'\i}nez Bernal for the 
painstaking discussion on the linear algebra related with this subject.
\newpage
\noindent
{\bsbf I. Basic properties of the Cayley--Dickson Algebras and its zero 
divisors}
\vglue.5cm
$\a_{n}$  denotes $\erre^{2^n}$ with the Cayley--Dickson multiplication 
$\,(n\ge 1)$:
\begin{eqnarray*}
x&=&(x_{1},x_{2}), \,y=(y_{1},y_{2})\quad\hbox{\bs in }\quad \erre^{2^n}=
\erre^{2^{n-1}}
\times \erre^{2^{n-1}} \\ 
xy&=&(x_{1}y_{1}-\overline{y}_{2}x_{2}, 
y_{2}x_{1}+x_{2}\overline{y}_{1})\quad \hbox{\bs with}\\ 
\overline{x}&=&(\overline{x}_{1},-x_{2})
\end{eqnarray*} 
$e_{0}=(1,0,\ldots,0)\in \a_{n}$ denotes the unit element. 
\par 
The Euclidian norm and inner product  are given by 
$$\parallel x\parallel^{2}=x\overline{x}=\overline{x}x=
\parallel x_{1}\parallel^{2}+\parallel x_{2}\parallel^{2}$$ 
and
$$\langle x,y\rangle={ 1\over 2 }(x\overline{y}+y\overline{x}),$$ 
respectively. 
\par 
\underline{The trace} is 
$$\hbox{\bs For}\,\, x\in \a_{n} \qquad t(x)=x+\overline{x},$$ 
i.e.
$$t(x)=2 (\hbox{\bs real part of} x)\quad \hbox{\bs and}\quad 
t(x)=2\langle x,e_{0}\rangle.$$
\par 
Let $x$ and $y$ be in $\a_{n}, x$ is orthogonal to 
$y\quad (x\perp y)\quad$ if and only
if \\
$x\overline{y}=-y\overline{x}\,$ or $\,t(x\overline{y})=0\,$ because
$\,\overline{xy} = \overline{y}\,\,\overline{x}$.
Therefore elements of zero trace (purely imaginaries) are orthogonal if and 
only if they anti--commute i.e. 
$$x\perp y\Longleftrightarrow xy=-yx$$ 
for
$x$ and $y$ in $\a_{n}$ with $t(x)=0$ i.e. $\overline{x}=-x$. 
\par 
Thus for all
$x\in \a_{n}$ we have the characteristic equation:
$$x^{2}-t(x)x+\parallel x\parallel^{2}=0$$
because
$$x^{2}-(x+\overline{x})x+\overline{x}x=0.$$
{\bsbf Definition:} The \underline{associator} of $a,b$ and $c$ in $\a_{n}$ is
$$(a,b,c)=(ab)c-a(bc)$$
which is linear in each variable.
\par 
For $n=0,1,2$ the associatior is indentically zero.
\par 
For $\a_{3}=\o$ the octonian numbers the associatior is alternative
i.e. for all $x$ and $y$ in $\a_{3}$%
$$x^{2}y=x(xy)\quad\hbox{\bs  and }\quad xy^{2}=(xy)y$$ 
then
\begin{eqnarray*}
0=(x+y,x+y,z)&=&(x,x,z)+(y,y,z)+(x,y,z)+(y,x,z)\, \hbox{\bs then}\\
0=(x,y,z)+(y,x,z)&\triangulo &(x,y,z)=-(y,x,z).
\end{eqnarray*}
Similarly $(x,y,z)=-(x,z,y)=(z,x,y)=-(z,y,x)$.
\par 
$\a_{n}$ for $n\ge 4$ is flexible i.e
$x(yx)=(xy)x$ for all $x$ and $y$ in $\a_{n}$. This means that 
$$(x,y,x)=0\quad\hbox{\bs or equivalently}\quad (x,y,z)=-(z,y,x)$$
for all $x, y$ and $z$ in $\a_{n}$.
\par 
Notice that $(x,y,z)=0$ if one of the entries 
$x,y$ or $z$ is real, that is, its imaginary part is zero.
\vglue.5cm
\noindent
{\bsbf Lemma 1.1.} For all $x,y$ and $z$ in $\a_{n}$

\begin{itemize}
\item[i)]$-(x,y,z)=(\overline{x},y,z)=(x,\overline{y},z)=(x,y,\overline{z})$
\item[ii)]$t((x,y,z))=0$
\item[iii)]$t(xy-yx)=0$
\end{itemize}
{\bsbf Proof.} For $n\le 2$ the assertions are trivial. Let $n\ge 3$ and 
$x=x_{r}+x_{I}$ where $x_{r}=$ real part of $x\quad$ and 
$x_{I}=$ imaginary part of
$x\quad$ so $\overline{x}=x_{r}-x_{I}$.
\begin{eqnarray*}
(\overline{x},y,z)&=&(x_{r}-x_{I},y,z)=(x_{r},y,z)-(x_{I},y,z)\\
&=&0-(x_{I},y,z)=-(x_{r},y,z)-(x_{I},y,z)\\
&=&-(x_{r}+x_{I},y,z)=-(x,y,z).
\end{eqnarray*}
Similarly $-(x,y,z)=(x,\overline{y},z)=(x,y,\overline{z})$. So we prove i).
\par 
To prove ii) we observe that
\begin{eqnarray*}
\overline{(x,y,z)}&=&\overline{-x(yz)}+\overline{(xy)z}=-
(\overline{z}\,\,\overline{y})\overline{x}+\overline{z}(\overline{y}\,\,
\overline{x})=\\
&=&-(\overline{z}, \overline{y}, \overline{x})=(z,y,x)=-(x,y,z)
\end{eqnarray*}
by i) and flexibility.
\par 
Therefore $(x,y,z)+\overline{(x,y,z)}=0$ so ii) is proved. 
\par 
If $x$ or $y$ are real
$xy=yx$ and $t(xy-yx)=0$. Suppose $x$ and $y$ are pure imaginary then
\begin{eqnarray*}
\overline{xy-yx}&=&\overline{xy}-\overline{yx}=\overline{y}\,\,\overline{x}-
\overline{x}\overline{y}=yx-xy\quad \triangulo \\
(xy-yx)&+&\overline{(xy-yx)}=xy-yx+yx-xy=0.
\end{eqnarray*}
On the other hand the symbol $[x,y]=xy-yx$ is bilinear, therefore decomposing
$x$ and $y$ in the real and imaginary parts respectively, we are done with 
(iii).
\begin{flushright}
Q.E.D.
\end{flushright}
\vglue.5cm
\noindent
{\bsbf Lemma 1.2.} For all $x,y,z$ and $w$ in $\a_{n}$.
$$x(y,z,w)+(x,y,z)w=(xy,z,w)-(x,yz,w)+(x,y,zw)$$
{\bsbf Proof.} (Adem's paper [1]).
\vglue.5cm
\noindent
{\bsbf Lemma 1.3.} For all $x,y$ and $z$ in $\a_{n}$%
$$\langle x,yz\rangle=\langle x\overline{z},y\rangle=
\langle\overline{y}x,z\rangle$$
{\bsbf Proof.}
\begin{eqnarray*}
2\langle x,yz\rangle&=&t(x(\overline{yz}))=t(x(\overline{z}\,\,\overline{y}))=
t((x\overline{z})\overline{y})\\
&=&2\langle x\overline{z},y\rangle
\end{eqnarray*}
by Lemma 1.1 ii).
\par 
Similarly
\begin{eqnarray*}
2\langle\overline{y}x,z\rangle&=&t((\overline{y}x)\overline{z})=t(
\overline{y}(x\overline{z}))=t((x\overline{z})\overline{y})=\\
&=&2\langle x\overline{z},y\rangle
\end{eqnarray*}
using Lemma 1.1. iii) and ii).
\begin{flushright}
Q.E.D.
\end{flushright}
\vglue.5cm
\noindent
{\bsbf Lemma 1.4.} For all $x$ and $y$ in $\a_{n}$.
$$\parallel xy\parallel = \parallel\overline{x}y\parallel$$
{\bsbf Proof.} 
$\parallel xy\parallel^{2}=\langle xy,xy\rangle$ and
$\parallel\overline{x}y\parallel^{2}=\langle\overline{x}y,\overline{x}y\rangle$
using Lemma 1.3.
\begin{eqnarray*}
\langle xy,xy\rangle&=&\langle\overline{x}(xy),y\rangle=
\langle(\overline{x}(xy))\overline{y},e_{0}\rangle={ 1\over 2 }t((\overline{x}(xy))
\overline{y})\\
\langle\overline{x}y,\overline{x}y\rangle&=&\langle x(\overline{x}y),y\rangle=
\langle(x(\overline{x}y))\overline{y},e_{0}\rangle=
{ 1\over 2 }t((x(\overline{x}y))
\overline{y}))
\end{eqnarray*}
thus 
\begin{eqnarray*}
\parallel xy\parallel^{2}-\parallel\overline{x}y\parallel^{2}&=&
{ 1\over 2 }t[(\overline{x}(xy)-x(\overline{x}y))\overline{y}]\\
&=&{ 1\over 2 }t[(-(\overline{x},x,y)+(\overline{x}x)y+(x,\overline{x},y)-
(x\overline{x})y)\overline{y}]\\
&=&{ 1\over  2}t[((x,x,y)+(x,\overline{x},y)+\parallel x\parallel^{2}y-
\parallel x\parallel^{2}y)\overline{y}]\\
&=&{ 1\over 2 }t(0\overline{y})=0.\\
&\triangulo &\parallel xy\parallel^{2}=\parallel\overline{x}y\parallel^{2}\quad
\hbox{\bs and }\quad \parallel xy\parallel=\parallel\overline{x}y\parallel
\end{eqnarray*}
\begin{flushright}
Q.E.D.
\end{flushright}
{\bsbf Corollary 1.5.} For all $x$ and $y$ in $\a_{n}$.
$$\parallel xy\parallel=\parallel\overline{x}y\parallel=
\parallel x\overline{y}\parallel=\parallel yx\parallel.$$
{\bsbf Proof.} Follows from the lemma recalling that 
$$\parallel z\parallel = \parallel\overline{z}\parallel\quad\forall 
z\in \a_{n}\quad
\hbox{\bs so}\quad
\parallel xy\parallel=\parallel\overline{x}y\parallel
=\parallel \overline{\overline{x}y}
\parallel=\parallel \overline{y}x\parallel=
\parallel yx\parallel$$
\begin{flushright}
Q.E.D.
\end{flushright}
\vglue.5cm
\noindent
{\bsbf Corollary 1.6.} (Elementary facts on zero divisors).
\par 
For $x$ and $y$ in $\a_{n}\quad n\ge 4$.
\par 
\begin{enumerate}
\item $xy=0 \Leftrightarrow yx=0\Leftrightarrow\overline{x}y=0
\Leftrightarrow x\overline{y}=0$.
\item If $x\neq 0$ and $xy=0$ then $t(y)=0$.
\item If $y\neq 0$ and $xy=0$ then $t(x)=0$.
\item $x^{2}=0$ if and only if $x=0$.
\end{enumerate}
{\bsbf Proof. } 1. By Corollary 1.5. 
$\parallel xy\parallel=\parallel yx\parallel=\parallel\overline{x}y\parallel
=\parallel x\overline{y}\parallel$ so
$xy=0\Leftrightarrow yx=0\Leftrightarrow\overline{x}y=0\Leftrightarrow x
\overline{y}=0$.
\par 
2. $xt(y)=x(y+\overline{y})=xy+x\overline{y}=0$\quad if $xy=0$\quad but $t(y)$
is a real number so $t(y)=0$.
\par 
3. $t(x)y=(x+\overline{x})y=xy+\overline{x}y=0\quad$ if $xy=0\quad$ 
but $t(x)$ is
a real number so $t(x)=0$.
\par 
4. If $x^{2}=x\cdot x=0$ then $t(x)=0\quad$ and $x=-\overline{x}\quad$ and 
$$\parallel x\parallel^{2}=\overline{x}x=-x^{2}\quad\triangulo\quad 
x^{2}=0\Leftrightarrow x=0.$$
\begin{flushright}
Q.E.D.
\end{flushright}
\vglue.5cm
\noindent
{\bsbf Notation:} $\circ{\a}_{n}=\{x\in \a_{n}|t(x)=0\}$.
\vglue.5cm
\noindent
{\bsbf Definition.} $0\neq x\in \circ{\a}_{n}$ is a zero divisor if there exists 
$0\neq y \in \circ{\a}_{n}$ with $xy=0$.
\par 
Notice that we don't have to distinguish between left and right zero divisors 
because $xy=0$ if and only if $yx=0$.
\par 
For $x\in \a_{n}$ we define the following linear transformations.
\begin{eqnarray*}
L_{x},R_{x}:\a_{n}&\rightarrow &\a_{n}\quad\hbox{by}\\
L_{x}(y)=xy&\hbox{\bs and}&R_{x}(y)=yx\quad\forall\quad y\in \a_{n}
\end{eqnarray*}
clearly $L_{x}$ and $R_{x}$ are linear.
\vglue.5cm
\noindent
{\bsbf Proposition 1.7.} For $x\in \circ{\a}_{n}\quad L_{x}$ and $R_{x}$ 
are skew--symmetric.
\vglue.5cm
\noindent
{\bsbf Proof.} Let $y$ and $z$ in $\a_{n}$. Then
\begin{eqnarray*}
\langle L_{x}(y), z\rangle&=&\langle xy,z\rangle=\langle y,\overline{x}z\rangle=
\langle y,-xz\rangle\\
&=&\langle y,-L_{x}(z)\rangle.\\
\langle R_{x}(y),z\rangle&=&\langle yx,z\rangle =\langle y,z\overline{x}\rangle=
\langle y,-zx\rangle\\
&=&\langle y,-R_{x}(z)\rangle.
\end{eqnarray*}
\begin{flushright}
Q.E.D.
\end{flushright}
\vglue.5cm
\noindent
{\bsbf Remark:} $x\in \circ{\a}_{n}$ is  a zero divisor if and only if
$\{0\}\neq $ Ker $L_{x}$ and \\
$\{0\}\neq$ Ker $R_{x}$ also notice that 
Ker $L_{x}=$ Ker $R_{x}$.
\vglue.5cm
\noindent
{\bsbf Lemma 1.8.} If $x\in \a_{n}$ is zero divisor and $x=(x_{1},x_{2})$ in
$\a_{n-1}\times \a_{n-1}$ then 
$
\check{x}\stackrel{\cdot\cdot}{=}(x_{2},x_{1})$ is zero divisor in 
$\a_{n}$. 
\vglue.5cm
\noindent
{\bsbf Proof.} Consider $y\neq 0$ in $\a_{n}$ such that $xy=0$ and 
$y=(y_{1},y_{2})$ so
$$xy=(x_{1},x_{2})(y_{1},y_{2})=(x_{1}y_{1}-\overline{y}_{2}x_{2},y_{2}x_{1}+x_{2}
\overline{y}_{1})=(0,0).$$
Consider the  following product in $\a_{n-1}\times \a_{n-1}=\a_{n}$%
\begin{eqnarray*}
(-\overline{y}_{2},y_{1})(x_{2},x_{1})&=&(-\overline{y}_{2}x_{2}-\overline{x}_{1}y_{1},
-x_{1}\overline{y}_{2}+y_{1}\overline{x}_{2})\\
&=&(x_{1}y_{1}-\overline{y}_{2}x_{2},\overline{x_{2}\overline{y}_{1}+y_{2}x_{1}})
\end{eqnarray*}
Recall that $t(x)=t(x_{1})=0$. Therefore $(-\overline{y}_{2},y_{1})(x_{2},x_{1})=(0,0)$
and \\
$\check{x}=(x_{2},x_{1})$ is a zero divisor.
\begin{flushright}
Q.E.D.
\end{flushright}
\vglue.5cm
\noindent
{\bsbf Corollary 1.9.} If $x\in \a_{n}$ is a zero divisor and $x=(x_{1},x_{2})$ then\\
$t(x)=t(x_{1})=t(x_{2})=0$.
\vglue.5cm
\noindent
{\bsbf Proof.} We already know that $t(x)=t(x_{1})=0$.
\par 
Now $\check{x}$ is a zero divisor then
$t(\check{x})=t(x_{2})=0$.
\begin{flushright}
Q.E.D.
\end{flushright}
\vglue.5cm
\noindent
{\bsbf Definition.} $a\in \a_{n}$ with $a=(a_{1},a_{2})\in \a_{n-1}\times \a_{n-1}$ is
{\bsit doubly pure} if $t(a_{1})=0$ and $t(a_{2})=0$.
\par 
So any zero divisor is doubly pure.
\par 
Let $e_{0}$ be the neutral element in $\a_{n-1}$ so the neutral element in 
$\a_{n}$ is $e_{0}=(e_{0},0)\in \a_{n-1}\times \a_{n-1}$.
\par 
Denote by $\widetilde{e}_{0}$ the element 
$\widetilde{e}_{0}\stackrel{\cdot\cdot}{=} (0,e_{0})\in \a_{n-1}\times \a_{n-1}$ i.e.
$\widetilde{e}_{0}=e_{2^{n-1}}$ in the canonical basis of $\erre^{2^n}$.
\vglue.5cm
\noindent
{\bsbf Lemma 1.10} Let $a\in \circ{\a}_{n}$ i.e. 
$t(a)=0$. Then $a$ is doubly pure if and
only if $\,a\perp \widetilde{e}_{0}$.
\vglue.5cm
\noindent
{\bsbf Proof.} Calculating 
$a\widetilde{e}_{0}=(a_{1},a_{2})(0,e_{0})=(-a_{2},a_{1})\,$ and\\
$\widetilde{e}_{0}a=(0,e_{0})(a_{1},a_{2})=(-\overline{a}_{2},\overline{a}_{1})=
(-\overline{a}_{2},-a_{1})\,$
because $\,t(a)=0$ and \\
$\overline{a}_{1}=-a_{1}$.
\par 
If $a$ is doubly pure $t(a_{2})=0$ and 
$\overline{a}_{2}=-a_{2}$ so 
$\widetilde{e}_{0}a=(a_{2},-a_{1})=-a\widetilde{e}_{0}$. Conversely if
$a\perp\widetilde{e}_{0}$ then $a\widetilde{e}_{0}=-\widetilde{e}_{0}a$ and
$(-a_{2},a_{1})=(\overline{a}_{2},a_{1})$ then
$\overline{a}_{2}=-a_{2}$ and $t(a_{2})=0$.
\begin{flushright}
Q.E.D.
\end{flushright}        
\vglue.5cm
Right multiplication by ``$\widetilde{e}_{0}$''  defines an orthogonal 
transformation of determinant one whose square is minus the identity linear 
transformation in $\a_{n}$.
\par 
For $x\in \a_{n}$ denotes by $\widetilde{x}=R_{\widetilde{e}_0}(x)=
(-x_{2},x_{1})$ if $x=(x_{1},x_{2})$ so
$$\parallel\widetilde{x}\parallel=\parallel x\parallel$$
and $\widetilde{\widetilde{x}}=-x$. 
\par 
And                      
$R^{T}_{\widetilde{e}_0}=-R_{\widetilde{e}_0}=R^{-1}_{\widetilde{e}_0}$ and 
$det(R_{\widetilde{e}_0})=(-1)^{2^n}=1$.
\vglue.5cm
\noindent
{\bsbf Proposition 1.11} For $a=(a_{1},a_{2})$ and
$b=(b_{1},b_{2})$ doubly pure elements in $\a_{n}$ we have:
\vglue.1cm
\begin{enumerate}
\item[1)]$a\widetilde{e}_{0}=\widetilde{a};\quad\widetilde{e}_{0}a=-\widetilde{a}$.
\item[2)]$a\widetilde{a}=-\parallel a\parallel^{2}\widetilde{e}_{0}$;\quad
$\widetilde{a}a=\parallel a\parallel^{2}\widetilde{e}_{0}\quad$ so
$\quad a\perp\widetilde{a}$.
\item[3)]$\widetilde{a}b=-\widetilde{ab}$.
\item[4)]If $a\perp b$ then $\widetilde{a}b+\widetilde{b}a=0$.
\item[5)]If $\widetilde{a}\perp b$ then $ab=\widetilde{b}\widetilde{a}$.
\end{enumerate}
{\bsbf Proof.} 1) Follows from definition.
\par 
2) $a\widetilde{a}=(a_{1},a_{2})(-a_{2},a_{1})
=(-a_{1}a_{2}+a_{1}a_{2},a^{2}_{1}+a^{2}_{2})
=(0,-\parallel a\parallel^{2})\\
=-\parallel a\parallel^{2}\widetilde{e}_{0}\\
\widetilde{a}a=(-a_{2},a_{1})(a_{1},a_{2})=(-a_{2}a_{1}+a_{2}a_{1},-a^{2}_{2}-a^{2}_{1})=
(0,\parallel a\parallel^{2})\\
=\parallel a\parallel^{2}\widetilde{e}_{0}$.
\par 
Since $t(a)=t(\widetilde{a})=0$ then
$a\widetilde{a}=-\widetilde{a}a$ and $a\perp\widetilde{a}$.
\par 
3) $\widetilde{a}b=(-a_{2},a_{1})(b_{1},b_{2})=
(-a_{2}b_{1}+b_{2}a_{1},-b_{2}a_{2}-a_{1}b_{1})\,$
so\\
$\,\widetilde{\widetilde{a}b}=(a_{1}b_{1}+b_{2}a_{2},b_{2}a_{1}-a_{2}b_{1})=ab
\quad\triangulo\quad \widetilde{a}b=-\widetilde{ab}$.
\par 
4) If $ab=-ba$ then $\widetilde{a}b=-\widetilde{ab}=\widetilde{ba}=-
\widetilde{b}a\quad\triangulo\quad \widetilde{a}b+\widetilde{b}a=0$.
\par 
5) If $\widetilde{a}\perp b\,$ then $0=\widetilde{ab}+
b\widetilde{a}=-\widetilde{ab}+\widetilde{ba}=ab+\widetilde{b\widetilde{a}}
=ab-\widetilde{ba}$ 
and
$\,0=\widetilde{\widetilde{a}}b+\widetilde{b}\widetilde{a}=-ab+\widetilde{b}
\widetilde{a}$.
\begin{flushright}
Q.E.D.
\end{flushright}
\vglue.5cm
\noindent
{\bsbf Corollary 1.12} For $0\neq a\in \a_{n}$ (doubly pure) \\
$ax=0\,$ if\, and only
if $\,a\widetilde{x}=0\quad\forall x$.
\vglue.5cm
\noindent
{\bsbf Proof.} $ax=0\Leftrightarrow xa=0\Leftrightarrow \widetilde{xa}
=0 \Leftrightarrow \widetilde{x}a=0\Leftrightarrow a\widetilde{x}=0$.
\begin{flushright}
Q.E.D.
\end{flushright}
\vglue.5cm
\noindent
{\bsbf Theorem 1.13} For $a\in \a_{n}$ doubly pure of norm one and $n\ge 2$.
The vector space generated by $\{e_{0},\widetilde{a},a,\widetilde{e}_{0}\}$ is
multiplicatively closed and isomorphic as algebra
to $\a_{2}=\ache$ the
quaternions.
\vglue.5cm
\noindent
{\bsbf Proof.} Construct the following multiplication table
\vglue.1cm
\begin{center}
\begin{tabular}{l|llll}
&$e_{0}$&$\widetilde{a}$&$a$&$\widetilde{e}_{0}$\\ \hline
$e_{0}$&$e_{0}$&$\widetilde{a}$&$a$&$\widetilde{e}_{0}$\\
$\widetilde{a}$&$\widetilde{a}$&$-e_{0}$&$\widetilde{e}_{0}$&$-a$\\              
$a$&$a$&$-\widetilde{e}_{0}$&$-e_{0}$&$\widetilde{a}$\\
$\widetilde{e}_{0}$&$\widetilde{e}_{0}$&$a$&$-\widetilde{a}$&$-e_{0}$
\end{tabular}
since
$
\begin{array}{lll}
a\widetilde{e}_{0}&=&\widetilde{a}\\
\widetilde{e}_{0}a&=&-\widetilde{a}\\
\widetilde{a}\widetilde{e}_{0}&=&\widetilde{\widetilde{a}}=-a\\
\widetilde{e}_{0}\widetilde{a}&=&-\widetilde{\widetilde{a}}=a\\
a\widetilde{a}&=&-\widetilde{e_{0}}\quad\hbox{\bs and}\quad
\widetilde{a}a=\widetilde{e_{0}}
\end{array}
$
\end{center}
\begin{flushright}
Q.E.D.
\end{flushright}
\vglue.5cm
\noindent
Recall that $Ker L_{ra}= Ker L_{a}$ for all $r\neq 0$ in $\erre$ and that 
$Ker L_{a}\neq \{0\}$ only if $a$ is doubly pure non--zero element for
$n\ge 4$.
\par 
{\bsit From now on we assume that $n\ge 4, \parallel a\parallel=1$ and
$a$ is doubly pure element in $\a_{n}$.}
\par 
$\ache_{a}$ denotes the copy of the quaternions inside of $\a_{n}$ generated
by $\{e_{0},\widetilde{a},a,\widetilde{e}_{0}\}$.
\par 
Let $T_{a}:\a_{n}\rightarrow \a_{n}$ be given by $T_{a}(y)=(a,a,y)$.
\par 
Because the associator is tri--linear, $T_{a}$ is a linear map and because 
$\ache_{a}$ is associative $T_{a}(\ache_{a})=\{0\}$.
\par 
On the other hand $T_{a}(y)=a^{2}y-a(ay)=-y-L^{2}_{a}(y)$ so\\
$T_{a}=-I-L^{2}_{a}$ and since $L_{a}$ 
is anti--symmetric then $T_{a}$ is symmetric \\
$\ache^{\perp}_{a}$ denotes the orthogonal complement of $\ache_{a}$ in
$\a_{n}=\erre^{2^n}$.
\vglue.5cm
\noindent
{\bsbf Lemma 1.14} $\quad T_{a}(\ache^{\perp}_{a})\subset \ache^{\perp}_{a}$.
\vglue.5cm
\noindent
{\bsbf Proof.} Let $0\neq y\in \ache_{a}^{\perp}\,$ then $\,T_{a}(y)=(a,a,y)$.
\par 
By Lemma 1.2.
$$a(a,a,y)+(a,a,a)y=(a^{2},a,y)-(a,a^{2},y)+(a,a,ay)\quad\triangulo$$
then $a(a,a,y)=(a,a,ay)$ because $a^{2}=-e_{0}$ and
$$\langle a,(a,a,y)\rangle =t(a\overline{(a,a,y)})=-t(a(a,a,y))=-t((a,a,ay))=0$$
because the trace of any asociator is zero (Lemma 1.1).
\par 
Similarly 
$\widetilde{a}(a,a,y)+(\widetilde{a},a,a)y=$ (Sum of Associators)  \\
since
$(\widetilde{a},a,a)=0\,$ because $\,\ache_{a}$ is associative then
$\widetilde{a}\perp (a,a,y)$.
\par 
Finally $\widetilde{e}_{0}(a,a,y)+(\widetilde{e}_{0},a,a)y=$ (Sum of Associators) 
since \\
$(\widetilde{e}_{0},a,a)=0\,$ because $\,\ache_{a}$ is associative then
$\widetilde{e}_{0}\perp (a,a,y)$.
\par 
Therefore $(a,a,y)\in \ache^{\perp}_{a}\,$ if $\,y\in \ache^{\perp}_{a}$
\begin{flushright}
Q.E.D.
\end{flushright}
\vglue.5cm
\noindent
$\widetilde{T}_{a}$ denotes the restriction of $T_{a}$ to $\ache_{a}^{\perp}\,$ so\\
$\widetilde{T}_{a}:\ache^{\perp}_{a}\rightarrow \ache^{\perp}_{a},\quad
\widetilde{T}_{a}=-I-L^{2}_{a}$ 
is symmetric.
\vglue.5cm
\noindent
{\bsbf Theorem 1.15} For each $a\in \a_{n}$ of norm one and doubly pure we have a 
direct sum decomposition.
$$\a_{n}=\erre^{2^n}\cong\ache_{a}\oplus Ker\widetilde{T}_{a}\oplus Ker L_{a}
\oplus \bigoplus V_{\lambda}$$
where
$$V_{\lambda}=\{x\in \a_{n}|a(ax)=-\lambda^{2}x\},
\quad\lambda>0,\lambda\neq 1$$
\vglue.5cm
\noindent
{\bsbf Proof.} Recall that $L_{a}$ is skew--symmetric and consequently $L^{2}_{a}$
is symmetric, negative semidefinite, that is for, all $x\in \a_{n}$%
$$\langle L^{2}_{a}(x),x\rangle=-\langle L_{a}(x), L_{a}(x)\rangle=-
|L_{a}(x)|^{2}\le 0$$
Therefore $a(ax)=0 \Leftrightarrow ax=0$ for $x\in \a_{n}$

On the other hand being $L^{2}_{a}$ symmetric no positive implies that all its 
eigenvalues are negative real numbers. So there exist 
$\{\lambda_{0},\lambda_{1},\ldots,\lambda_{s}\}\subset \erre$ with 
$\lambda_{0}=0$ and $\lambda_{1}=1$ and eigenspaces
$$V_{\lambda}=\{x\in \ache^{\perp}_{a}|a(ax)=-\lambda^{2}x\}.$$
Recall that $L^{2}_{a}|\ache_{a}=-I$.
\par 
And $V_{0}= Ker L^{2}_{a}= Ker L_{a}$ and $V_{1}= Ker \widetilde{T}_{a}$.
\begin{flushright}
Q.E.D.
\end{flushright}
\vglue.5cm
\noindent
{\bsbf Theorem 1.16.} For the $V_{\lambda}$'s as in the last theorem we have
that for \\ $\lambda>0$%
$$ dim_{\erre}V_{\lambda}\equiv 0\quad mod\,\, 4.$$
\vglue.5cm
\noindent
{\bsbf Proof.} 
\par 
{\bsbf CLAIM 1.} For $x\in V_{\lambda}$ with 
$\lambda >0\quad y\stackrel{\cdot\cdot}{=}\lambda^{-1}(ax)$  
belongs $V_{\lambda}$, and $y\perp x$.
\par 
$ay=\lambda^{-1} a(ax)=\lambda^{-1}(-\lambda^{2}x)=-\lambda x$ and
$a(ay)=-\lambda ax=-\lambda^{2}(\lambda^{-1} ax)=-\lambda^{2}y$.
Now $\langle y,x\rangle=\langle\lambda^{-1} ax,x\rangle=\lambda^{-1}\langle 
ax,x\rangle=0$.
\par 
{\bsbf CLAIM 2.} If $x\in V_{\lambda}$ then $\widetilde{x}\in V_{\lambda}$.
\par 
For this we observe that: Since 
$x\in \ache^{\perp}_{a}, a\perp x, a\perp\widetilde{x}$ and $x$ is doubly pure then 
applying Proposition 1.11 (3) twice and noting that $a(ax)=(xa)a$ 
by flexibility we have
that 
$a(a\widetilde{x})=(\widetilde{x}a)a=-(\widetilde{(xa)a})=\widetilde{a(ax)}=
-\lambda^{2}\widetilde{x}$.
\par 
Therefore if we take $0\neq x\in V_{\lambda}$ for $\lambda\neq 0$ we may construct 
one orthogonal set
$$\{x,y,\widetilde{x},\widetilde{y}\}$$
inside of $V_{\lambda}$. Thus $dim_{\erre}V_{\lambda}\equiv 0\quad
mod\,\, 4 \quad \lambda >0$
\begin{flushright}
Q.E.D.
\end{flushright}
\vglue.5cm
\noindent
{\bsbf Corollary 1.17.}
\par 
$$2^{n}-4\ge dim_{\erre} Ker L_{a}\equiv 0\quad mod\,\, 4$$
\vglue.5cm
\noindent
{\bsbf Proof.} Since $dim \ache_{a}=4, 2^{n}-4\ge  dim_{\erre}Ker L_{a}$ and
the last two theorems implies that
$$dim_{\erre}Ker L_{a}\equiv 0\quad mod\,\, 4.$$
\begin{flushright}
Q.E.D.
\end{flushright}
\vglue.5cm
\noindent
\newpage
\vglue.5cm
\noindent
{\bsbf II. Zero divisors with alternative entries}
\vglue.5cm
Throughout this chapter $a$ and $b$ are nonzero elements in $\a_{n}$ for
$n\ge 3$.
\vglue.5cm
\noindent
{\bsbf Definition:} $a\in \a_{n}$ is \underline{alternative} if $(a,a,x)=0$ for all
$x\in \a_{n}$.
\par 
The real elements $(\erre e_{0})$ are alternative and since the associator symbol 
is linear in each variable we have that one element in $\a_{n}$ is alternative if 
and only if it's imaginarie part is alternative. 
\par 
Therefore we restrict 
ourselves to the trace zero of norm one alternative elements in $\a_{n}$.
\par 
By Schafer [8] we know that the canonical basis $\{e_{0},e_{1},\ldots,e_{2^n-1}\}$
consists of alternative elements.
\vglue.5cm
\noindent
{\bsbf Lemma 2.1.} If $a\in \a_{n}$ is alternative with $t(a)=0$ and $|a|=1$ then
\par 
\begin{enumerate}
\item[i)] $Ker L_{a}=\{0\}$
\item[ii)]$L_{a}\in SO(2^{n})$
\item[iii)]$ax=y$ has a unique solution for all $y\in \a_{n}$
\item[iv)]$L^{2}_{a}=-I$
\end{enumerate}
\vglue.5cm
\noindent
{\bsbf Proof.} $0=(a,a,x)=-x-a(ax)$ then $a(ax)=-x$ and if
$ax=0$ then $a(ax)=-x=0$ and $Ker L_{a}=\{0\}$ and $L_{a}$ is one to one.
Now $L_{a}$ is skew--symmetric so
\begin{eqnarray*}
0&=&\langle(a,a,x),x\rangle=\langle -x-a(ax),x\rangle=-\langle x,x\rangle-
\langle a(ax),x\rangle\\
&=&-|x|^{2}+\langle ax,ax\rangle=-|x|^{2}+|ax|^{2}
\end{eqnarray*}
then $|L_{a}(x)|=|x|$ for all $x\in \a_{n}$ and $L_{a}\in O(n)$. Therefore
$L_{a}^{-1}=L^{T}_{a}=-L_{a}$ and $L^{2}_{a}=-I$ and $L_{a}$ is similar to
$
\left( 
\begin{array}{cc}
0&I\\ -I&0
\end{array}\right)
$ 
for $I$ and $0$ in
$M_{2^{n-1}}(\erre)$ which has determinant one.
\begin{flushright}
Q.E.D.
\end{flushright}
\vglue.5cm 
\noindent
Also we have similar results for $R_{a}$, right multiplication.
\vglue.3cm
\noindent
{\bsbf Definition:} $a\in \a_{n}$ is \underline{special} if $a$ is alternative, $t(a)=0$
and $|a|=1$.
\par 
If $a\in \a_{n}$ is special then $\{e_{0},a\}$ generates a copy a $\a_{1}=\ce$ inside
of $\a_{n}$. So the special elements can be regarded as a closed subset of the 
manifold of the $2^{n}\times 2 $ matrices of rank two with real coeficients.
\par 
\vglue.5cm
\noindent
{\bsbf Definition.} $\{a,b\}\subset \a_{n}$ is a special couple if both $a$ and
$b$ are special elements in $\a_{n}$ and $a$ is orthogonal to $b\quad (a\perp b)$
i.e. $ab=-ba$.
\vglue.5cm
\noindent
{\bsbf Notation:} $V(a;b)$ denotes the vector space generated by the set
$$\{e_{0},a,b,ab\},$$
obviously $V(a;b)$ is fourth--dimensional: $a\perp ab$ and $b\perp ab$.
\vglue.5cm
\noindent
{\bsbf Proposition 2.2.} $V(a;b)$ 
is multiplicatively closed and isomorphic to\\
$\a_{2}=\ache$ the quaternions.
\vglue.5cm
\noindent
{\bsbf Proof:} Since $a\perp b, t(a\overline{b})=-t(ab)=0$ and
$|ab|^{2}=\langle ab,ab\rangle=\langle -a(ab),b\rangle=\langle b,b\rangle=
|b|^{2}=1$. 
Therefore $t(ab)=0$ and $|ab|=1$. $\{e_{0},a,b,ab\}$ is an orthonormal set.
Clearly $a(ab)=a^{2}b=-b$ and $b(ab)=-b(ba)=-b^{2}a=a$. So we have the following 
multiplication table
\vglue.1cm
\begin{center}
\begin{tabular}{l|rrrr}
&$e_{0}$&$a$&$b$&$ab$\\ \hline
e$_{0}$&$e_{0}$&$a$&$b$&$ab$\\
$a$&$e_{0}$&$-e_{0}$&$ab$&$-b$\\              
$b$&$b$&$-ab$&$-e_{0}$&$a$\\
$ab$&$ab$&$b$&$-a$&$-e_{0}$
\end{tabular}
\end{center}
\vglue.2cm
\begin{flushright}
Q.E.D.
\end{flushright}
\vglue.5cm
Therefore each special couple induce a multiplicative monomorphism from
$\a_{2}=\ache$ to $\a_{n}$ and the set of special couples can be regarded as 
closed subset of the $2^{n}\times 4$ matrices with real entries and rank four.
\par 
Depending on $\{a,b\}$ special couple we define
$$S:\a_{n}\rightarrow \a_{n}\quad\hbox{\bs by}\quad S(y)=(a,y,b)$$
First of all $S$ is a linear map and $V(a;b)\subset Ker S$ because
$V(a;b)\cong\ache$ is associative. Also $S=R_{b}L_{a}-L_{a}R_{b}=[R_{b},L_{a}]$ is
skew--symmetric because $R_{b}$ and $L_{a}$ are skew--symmetric.
\par 
Therefore $S^{2}(y)=S(S(y))=0$ if and only if $S(y)=0$ and
$S:V(a;b)^{\perp}\rightarrow V(a;b)^{\perp}$ because $S(y)\in V(a;b)$ for $y\in V(a;b)$.
We are interested in calculate $Ker S$ and $Im S$ for $\{a,b\}$ special couple
and $S$ restricted to $V(a,b)^{\perp}$.
\vglue.5cm
\noindent
{\bsbf Theorem 2.3.} For $\{a,b\}$ special couple and $S$ and $V(a,b)$ as above
we have that 
$$Ker(S:V(a;b)^{\perp}\rightarrow V(a;b)^{\perp})=Ker L_{a+b}\oplus Ker L_{a-b}.$$
{\bsbf Proof.} First of all we notice that $Ker L_{a+b}\cap Ker L_{a-b}=\{0\}$.
If $0=(a+b)z=(a-b)z$ then $az+bz=az-bz=0$ and $2bz=0$ and $z=0$ because
$b$ is alternative. Now if $(a+b)x=0$ then $ax=-bx=xb$ for $x\in V(a;b)^{\perp}$
then $(a,x,b)=(ax)b-a(xb)=(xb)b-a(ax)=xb^{2}-a^{2}x=0$. Similarly if
$(a-b)x=0$ then $ax=bx=-xb$ for $x\in V(a,b)^{\perp}$. Thus
$(a,x,b)=(ax)b-a(xb)=-(xb)b+a(ax)=-xb^{2}+a^{2}x=0$. Therefore
$Ker L_{a+b}\oplus Ker L_{a-b}\subset Ker S$.
Now if $y\in V(a;b)^{\perp}$ with $S(y)=0$ we have that
\begin{eqnarray*}
(a+b)(y-a(yb))&=&ay+by-a^{2}yb-b(a(yb))\\
&=&ay+by-by+((ay)b)b\\
&=&ay+ayb^{2}\\
&=&ay-ay\\
&=&0
\end{eqnarray*}
because $a(yb)=(ay)b$ and $b\perp[(ay)b]$.
\par 
Similarly $(a-b)(y+a(yb))=0$ if $S(y)=0$ and $y\in V(a;b)^{\perp}$. Since
$y={ 1\over  2}[(y-a(yb))+(y+a(yb)]$ we have that: If $S(y)=0$ and
$y\in  V(a;b)^{\perp}$ then $y\in Ker L_{a+b}\oplus Ker L_{a-b}$ and we are done.
\begin{flushright}
Q.E.D.
\end{flushright}
\vglue.5cm
\noindent
{\bsbf Corollary 2.4} $dim_{\erre} Ker S=2 dim Ker L_{a+b}\equiv 0$ mod 8.
\vglue.5cm
\noindent
{\bsbf Proof.} By the proof of the Theorem 2.3. we know that the elements in
$Ker L_{a+b}$ are of the form $[y-(a(yb)]$ and the elements in $Ker L_{a-b}$ are
of the form $[y+(a(yb))]$ for $(a,y,b)=0$. Therefore the assignment
$$y+a(yb)\mapsto y-a(yb)$$
defines a linear isomorphism between $Ker L_{a+b}$ and $Ker L_{a-b}$.
Therefore $dim_{\erre} Ker S=2 dim Ker L_{a+b}$ that by Theorem 1.15 is congruent
with $0$ module $8$.
\vglue.5cm
\noindent
{\bsbf Lemma 2.5} For $\{a,b\}$ special couple in $\a_{n}$ and $y\in V(a;b)^{\perp}$ 
we have that
$$(a+b,a+b,y)=-(a,y,b)+2(ay)b=(a,y,b)+2a(yb)$$
in short $T_{a+b}=-S+2R_{b}L_{a}=S+2L_{a}R_{b}$ on $V(a;b)^{\perp}$
\vglue.5cm
\noindent
{\bsbf Proof.}
\begin{eqnarray*}
(a+b,a+b,y)&=&(a,a,y)+(b,b,y)+(a,b,y)+(b,a,y)\\
&=&0+0+(ab+ba)y-a(by)-b(ay)\\
&=&a(yb)+(ay)b\\
&=&-(a,y,b)+2(ay)b\\
&=&(a,y,b)+2a(yb)
\end{eqnarray*}
\begin{flushright}
Q.E.D.
\end{flushright}
\vglue.5cm
\noindent
{\bsbf Corollary 2.6} Let $\{a,b\}$ be one special couple and $y\in V(a,b)^{\perp}$
in $\a_{n}$.
$(a+b,a+b,y)=0\quad$ if and only if
$$(a,b)(-{ 1\over 2 }S(y),y)=(0,0)\quad\hbox{\bs in}\quad \a_{n+1},$$
{\bsbf Proof.} By Lemma 2.5 
$T_{a+b}(y)=0\Leftrightarrow a(yb)=-(ay)b$ then
\begin{eqnarray*}
(a,b)(a(yb),y)&=&(a(a(yb))+yb,ya-b(a(yb)))\\
&=&(a^{2}yb+yb,ya-((ay)b)b)\\
&=&(-yb+yb,ya-(ay)b^{2})\\
&=&(0,0)
\end{eqnarray*}
because $y\in V(a,b)^{\perp}$.
\par 
Conversely if $(a,b)(x,y)=(0,0)$ for $y\in V(a;b)^{\perp}$ then
$ax+yb=0$ and $ya-bx=0$ so
\begin{eqnarray*}
ax=-yb&\Rightarrow& a(ax)=a^{2}x=-x=-a(yb)\,\,\triangulo\,\, x=a(yb)\\
bx=ya&\Rightarrow&ay=xb\Rightarrow(ay)b=(xb)b=xb^{2}=-x.
\end{eqnarray*}
Therefore $a(yb)+(ay)b=x-x=0$. And by the Lemma 2.5 $T_{a+b}(y)=0$ and
$x={ 1\over  2}S(y)$
\begin{flushright}
Q.E.D.
\end{flushright}
\par 
\vskip.2cm
\par 
\noindent
{\bsbf Theorem 2.7} Let $\{a,b\}$ be one special couple in $\a_{n}$ and \\
$L_{(a,b)}:\a_{n+1}\rightarrow \a_{n+1}\quad$ then
$$dim_{\erre} Ker L_{(a,b)}\le 2^{n}-4-2 dim_{\erre}Ker L_{a+b}$$
{\bsbf Proof. } By the last corollary
$$Ker L_{(a,b)}\cong Ker[T_{a+b}:V(a;b)^{\perp}\rightarrow V(a;b)^{\perp}]$$
but $Ker T_{a+b}\subset Im S$ restricted to $V(a;b)^{\perp}$ and
$Im S\cong(Ker S)^{\perp}=(Ker L_{a+b}\oplus Ker L_{a-b})^{\perp}$
(by Theorem 2.3).
Also we know that $Ker L_{a+b}\cong Ker L_{a-b}$ and
$dim_{\erre} Ker L_{(a,b)}\le dim_{\erre}\a_{n}-dim_{\erre}V(a,b)-2 dim_{\erre} Ker L_{a+b}$.
\begin{flushright}
Q.E.D.
\end{flushright}
\vglue.5cm
\noindent
{\bsbf Examples:}
\par 
For $n=3$. All element in $\a_{3}$ are alternative so 
$Ker T_{a+b}=\a_{3}$ for all
$a$ and $b$. Then $dim Ker L_{(a,b)}=4$ for all special couple with
$S(y)\neq 0$ for some $y\in V(a;b)^{\perp}$.
\par 
For $n\ge 4$. The top dimension $2^{n}-4$ is always realizable. For instance if $a$
is one special element and doubly pure in $\a_{n}$ then $\widetilde{a}$ is also
special and $V(a;\widetilde{a})=\ache_{a}$ then $Ker T_{a+\widetilde{a}}=\{0\}$
so $dim Ker L_{(a,\widetilde{a})}=2^{n}-4$.
\par 
More generally if $\{a,b\}$ is one special couple in $\a_{n}$ with $(a+b)$
alternative then $Ker L_{(a,b)}\cong V(a;b)^{\perp}$ i.e. $dim Ker L_{(a,b)}=2^{n}-4$.
\par 
Now we analyze a more general case:
\par 
Given $a$ and $b$ alternative non--zero elements in $\a_{n}$.
\par 
{\bsit 
Under what conditions on $a$ and $b, (a,b)$ is a zero divisor in $\a_{n+1}$?}
\par 
Suppose that $(a,b)(x,y)=(0,0)$ for $x\neq 0$ and $y\neq 0$. Then
\begin{enumerate}
\item[1)] $t(a)=t(b)=0$ (Corollary 1.9)
\item[2)] $|a|=|b|$\\
$ax=-yb\Rightarrow |a||x|=|b||y|$ and\\
$ya=bx\Rightarrow |a||y|=|x||b|$ so
$$|a|^{2}|y|=|a||x|b|=|b|^{2}|y|\quad\hbox{\bs and}\quad |a|=|b|.$$
Notice that without loosing generality we may assume that
$$|a|=|b|=1$$ 
from now on.
\item[3)] $a$ and $b$ are linearly independent.
\par 
Suppose that $a=\lambda b$ for $\lambda\neq 0$ in $\erre$ so\\
$ax=-yb\Rightarrow \lambda bx=-yb$ and\\
$ya=bx\Rightarrow \lambda yb=bx$ so
$$\lambda (bx)=\lambda(\lambda yb)=\lambda^{2}yb=-yb\Rightarrow \lambda^{2}=-1$$
which is a contradiction.
\item[4)] $(a,y,b)=-2x$ and $(a,x,b)=2y$. Recall that 
$(a,b)(x,y)=(0,0)\Leftrightarrow (a,b)(-y,x)=(0,0)$. Thus
\begin{eqnarray*}
ax=-yb&\hbox{\bs and}&xa=-by\\
ya=bx&&ay=xb
\end{eqnarray*}
Then
\begin{eqnarray*}
ax=-yb&\Rightarrow&a^{2}x=-a(yb)\Rightarrow x=a(yb)\\
ay=xb&\Rightarrow&(ay)b=xb^{2}\Rightarrow x=-(ay)b
\end{eqnarray*}
\end{enumerate}
Also
\begin{eqnarray*}
xa=-by&\Rightarrow&b(xa)=-b^{2}y=y\\
ya=bx&\Rightarrow&ya^{2}=(bx)a\Rightarrow -y=(bx)a
\end{eqnarray*}
Therefore $(a,y,b)=(ay)b-a(yb)=-x-x=-2x$ and
$(a,x,b)=-(b,x,a)=-[(bx)a-b(xa)]=2y$.
\par 
Notice that we also show that $x\perp y$ because 
$S\stackrel{\cdot\cdot}{=}(a,-,b)$ is skew--symmetric and that $x$ and $y$
belongs to $\{e_{0},a,b,ab\}^{\perp}$ because 
$S(e_{0})=S(a)=S(b)=0$ and $S(ab)=-aS(b)=0$.
\begin{enumerate}
\item[5)]$T_{a+b}(y)\stackrel{\cdot\cdot}{=}(a+b,a+b,y)=-2\langle a,b\rangle y$
and
\begin{eqnarray*}
L^{2}_{a+b}(y)&\stackrel{\cdot\cdot}{=}&(a+b)[(a+b)y]=-2y\\
(a+b)[(a+b)y]&=&(a+b)[ay+by]=a(ay)+b(by)+b(ay)\\
&&\qquad\qquad\qquad\qquad +a(by)\\
&=&-y-y+x-x=-2y
\end{eqnarray*}
because $b(ay)=-(ay)b=x$ and $a(by)=-a(yb)=-x$.
\end{enumerate}
And 
\begin{eqnarray*}
T_{a+b}(y)&=&(a+b)^{2}y-L^{2}_{a+b}(y)\\
&=&(a^{2}+b^{2}+ab+ba)y+2y\\
&=&-2y-2\langle a,b\rangle y+2y\\
&=&-2\langle a,b\rangle y
\end{eqnarray*}
We collect enough necessary conditions to stablish.
\vglue.5cm
\noindent
{\bsbf Theorem 2.9} Let $a$ and $b$ alternatives elements in $\a_{n}$ of
zero trace and norm one, then $(a,b)\in \a_{n+1}$ is a zero divisor if and only
if $\lambda=-2$ is an eigenvalue of $L^{2}_{a+b}$
\vglue.5cm
\noindent
{\bsbf Proof.}
\par 
$\Rightarrow$) Follows from the conditions enlisted above.
\par 
$\Leftarrow$) Let $y\neq 0$ in $\a_{n}$ such that $L^{2}_{a+b}(y)=-2y$ then
\begin{eqnarray*}
-2y&=&(a+b)[(a+b)y]=a^{2}y+by^{2}+a(by)+b(ay)\\
&=&-2y+a(by)+b(ay)
\end{eqnarray*}
therefore 
$$a(by)=-b(ay).$$
On the other hand
\begin{eqnarray*}
L^{2}_{a-b}(y)&=&(a-b)[(a-b)y]=(a^{2}+b^{2})y-(b(ay)+a(by))\\
&=&-2y.
\end{eqnarray*}
then $y\perp (a+b)$ and $y\perp(a-b)$ and $y\perp a$ and $y\perp b$

Consider the following product in $\a_{n+1}$
\begin{eqnarray*}
(a,b)(-a(by),y)&=&(-a(a(by))+yb,ya+b(a(by))\\
&=&(-a^{2}(by)+yb,ya-b(b(ay))\\
&=&(by+yb,ya-b^{2}(ay)\\
&=&(-yb+yb, ya-ya)\\
&=&(0,0)
\end{eqnarray*}
because $a(by)=-b(ay), a\perp y$ and $b\perp y$.
\begin{flushright}
Q.E.D.
\end{flushright}
\vglue.5cm
\noindent
{\bsbf Remark:} Condition 5) above shows that Theorem 2.9 implies Corollary 
2.6 it is the case when $\langle a,b\rangle=0$.
\vglue.5cm
\noindent
{\bsbf Corollary 2.10} For $a$ and $b$ alternative elements in $\a_{n}$ of
zero trace and norm one
$$dim Ker L_{(a,b)}\le 2^{n}-4-2 dim Ker L_{a+b}$$
\newpage
\noindent
{\bsbf Proof.} 
\par 
Notice that $a(by)=-b(ay)\Leftrightarrow S(y)=-2(ay)b$ for
$y\in V(a,b)^{\perp}$ such that $(a,b)(x,y)=(0,0)$ then
$$Ker L_{(a,b)}\subset (Ker S)^{\perp}\cong (Ker L_{a+b}\oplus Ker L_{a-b})^{\perp}$$
then $dim Ker L_{(a,b)}\le 2^{n}-4-2 dim Ker L_{a+b}$ because 
$Ker L_{a+b}\cong Ker L_{a-b}$.
\begin{flushright}
Q.E.D.
\end{flushright}
\vglue.5cm
\noindent
{\bsbf Remarks:} Notice $L^{2}_{a+b}(y)=-2y$ for $a,b$ as in the Theorem 2.9
and $y\neq 0$ implies that $a$ and $b$ are linearly independent. If it would
$a=b$ or $a=-b$ then $L^{2}_{a+b}=-4I$ and $L^{2}_{a-b}=0$.
\par 
There are examples of couples of alternatives 
linearly independent elements (of norm one and zero 
trace) in $\a_{n}$
which no form a zero divisor in $\a_{n+1}$. For instances put $a=e_{1}$
and $b={ e_{1}+e_{2}\over  \sqrt {2 }}$ in $\a_{3}$ then $(a,b)\in \a_{4}$ is neither a
zero divisor nor an alternative element. (We thank to Paul Yiu for this piece 
of information).
\vglue.5cm
\noindent
{\bsbf Definition.} A triple 
$\{a,y,b\}$ in $\a_{n}$ is special if it is an
orthonormal set and 
\par 
i) $(a,b)$ is an special couple
\par 
ii) $(ay)b =-a(yb)$
\vglue.5cm
\noindent
{\bsbf Definition.} A zero divisor $(a,b)$ in $\a_{n+1}$ is special if 
$(a,b)(x,y)=(0,0)$ implies that $\{a,y,b\}$ is one special triple.
\vglue.5cm
\noindent
{\bsbf Proposition 2.11} Let $a$ and $b$ are alternatives of norm one in $\a_{n}$
with 
$(a,b)(x,y)=(0,0)$ in $\a_{n+1}$. If $(a+b)$ is alternative then $(a,b)$
is an special zero divisor.
\vglue.5cm
\noindent
{\bsbf Proof.} By Theorem 2.9) $T_{a+b}(y)=-2\langle a,b\rangle y$. So
$0=T_{a+b}(y)\Rightarrow a\perp b$ and $\{a,b\}$ is one special couple on the other 
hand $(ay)b=-(ay)b$ by 4).
\begin{flushright}
Q.E.D.
\end{flushright}
\vglue.5cm
\noindent
{\bsbf Corollary 2.12} Any zero divisor (up to norm) in $\a_{4}$ is special
zero divisor.
\par 
Consider the vector subspace in $\a_{n}$%
$$V(a;y;b)\stackrel{\cdot\cdot}{=}\langle \{e_{0},a,b,ab,(ay)b,yb,ay,y\}\rangle$$
\vglue.3cm
\noindent
{\bsbf Theorem 2.13} $V(a;y,b)$ is a copy of $\a_{3}=\o$ the octonian numbers
inside of $\a_{n}$ if and only if $\{a,y,b\}$ is an special triple.
\vglue.3cm
\noindent
{\bsbf Proof.} Let's make the following assignment.
$$
\begin{array}{cccccccc}
e_{0}&a&b&ab&(ay)b&yb&ay&y\\[.5cm]
\rule{.01cm}{1cm}&
\rule{.01cm}{1cm}&
\rule{.01cm}{1cm}&
\rule{.01cm}{1cm}&
\rule{.01cm}{1cm}&
\rule{.01cm}{1cm}&
\rule{.01cm}{1cm}&
\rule{.01cm}{1cm}\\
e_{0},&e_{1},&e_{2},&e_{3},&e_{4} &e_{5} &e_{6} &e_{7}
\end{array}
$$
we may easily see that this define an homomorphism of $\a_{3}$ to $\a_{n}$ if
$\{a,y,b\}$ is an special triple. Notice in fact that this a monomorphism from
$\a_{3}$ to $\a_{n}$.
\par 
Conversely is easy to see that $\{e_{1},e_{7}, e_{2}\}$ form  in special triple in
$\a_{n} $ for $n\ge 3$.
\begin{flushright}
Q.E.D.
\end{flushright}
\par 
\vskip.5cm
\par 
\noindent
{\bsbf Corollary 2.14} The set of zero divisor in $\a_{4}$ (with entries of norm
one) are homemorphic to $Aut(\a_{3})=G_{2}$ the exceptional Lie group of rank 
two.
\vglue.5cm
\noindent
{\bsbf Remarks:} Notice that if $\{a,y,b\}$ is an special triple then 
$\{b,y,a\}$ is an special triple and $\{a,b,y\}$ is no necessarily an special
triple. e.g. \\
$a=e_{1}\quad b=e_{2}\quad y=e_{15}$ in $\a_{4}$.
$(e_{1},e_{15},e_{2})=2e_{12}$ and \\
$(e_{1}e_{15})e_{2}=-e_{1}(e_{15}e_{2})=e_{12}$
but $(e_{1},e_{2},e_{15})=0$ and $(e_{1}e_{2})e_{15}=e_{1}(e_{2}e_{15})$.
\vglue1cm
\noindent
{\bsbf REFERENCES}
\vglue.5cm
\baselineskip.3cm
\begin{enumerate}
\item Adem, Construction of Some Normed Maps. Bol. Soc. Mat. Mexicana 
(1975), 59--75.
\item F.R. Cohen. On Whitehead Squares, Cayley--Dickson algebras and 
rational functions. Bol. Soc. Mat. Mexicana 37(1992), 55--62.
\item P. Eakin, A. Sathaye. On Automorphisms and Derivations of Cayley--Dickson
Algebras. Journal of Algebra (1990).
\item P. Eakin, A. Sathaye, Yiu's Conjecture. Unpublished, University of 
Kentucky (1992).
\item Kantor--Solodovnikov. Hypercomplex Numbers. Springer--Verlag 
(1989).
\item S. Khalil. The Cayley--Dickson Algebras. (Master's Thesis 1993), 
Florida Atlantic University.
\item Koecher--Remmert. Numbers. Part B. Graduate Texts in Mathematics 
123, Springer--Verlag 1991.
\item R.D. Schafer. On algebras formed by the Cayley--Dickson process.
American Journal of Mathematics, (1954), 435--446.
\item W. Whitehead. Elements of Homotopy Theory. (Appendix A). Graduate
text in Mathematics 61, Springer--Verlag.
\end{enumerate}
\bf
Departamento de Matem\'aticas\\
Centro de Investigaci\'on y de Estudios Avanzados del I.P.N. \\
Apdo. Post. 14-740,07000, M\'exico, D.F., M\'exico \\
e-mail: gmoreno@math.cinvestav.mx 

\end{document}